%&amstex
\input amstex.tex
\documentstyle{amsppt}
\magnification=1200
\hsize=150truemm
\vsize=224.4truemm
\hoffset=4.8truemm
\voffset=12truemm

\TagsOnRight
\NoBlackBoxes
\NoRunningHeads

\refstyle{A}

\topmatter

\title On modular forms of characteristic $p>0$.
\endtitle

\author Marc Reversat \endauthor

\affil Laboratoire de Math\'ematiques \'Emile Picard, Unit\'e Mixte de Recherche
Universit\'e-C.N.R.S. 5580. 
\endaffil

\address {Universit\'e Paul Sabatier, 118 route de Narbonne, 31062 Toulouse c\'edex 4. France}
\endaddress

\email reversat\@picard.ups-tlse.fr \endemail

\date{}
\enddate

\abstract{We compare modular forms of characteristic $p>0$ (i.e. Drinfeld's modular forms) 
and automorphic forms. We prove that spaces of these modular forms (which are of characteristic
$p$) can be described by function spaces of characteristic zero, close to those of
automorphic forms.}
\endabstract

\endtopmatter

\document

\head 0 Introduction.\endhead

\noindent {\bf (0.1)} Let $K$ be a global field of characteristic $p>0$ (i.e. a function field
of one variable over a finite field of characteristic $p$) with a marked place, denoted by
$\infty$. For any place $v$ of $K$, we denote by $K_v$ the completion of
$K$ at $v$, by ${\Cal O}_v$, the valuation ring of $K_v$. 
Let $A$ be the subring of $K$ of regular elements away from $\infty$
(i.e. of $\lambda \in K$ such that $\lambda \in {\Cal O}_v$ for all $v\not= \infty$). 

\vskip5pt \noindent {\bf (0.2)} $G$ denotes the group-scheme $GL_2$ and $Z$ is its center.

\vskip5pt \noindent
{\bf (0.3)} The ring of ad\`eles of $K$, denoted by ${\Bbb A}$, can be written
${\Bbb A}={\Bbb A}_f\times K_\infty$, where ${\Bbb A}_f$ is the restricted product of 
$\{K_v\}_v$,
 $v$ running over the set of places of $K$ not equal to $\infty$ (the elements of ${\Bbb A}_f$
are called finite ad\`eles). One sets also
$\Cal O=\prod_v {\Cal O}_v$ ($v$ runs over the set of all places of $K$) and
${\Cal O}_f=\prod_{v\not= \infty}{\Cal O}_v$.

\vskip5pt \noindent
{\bf (0.4)} Following G. Harder, we will underline elements of adelic nature: for
instance, an element ${\underline g}\in G({\Bbb A})=G({\Bbb A}_f)\times G(K_\infty)$ may be
decomposed as  ${\underline g}=({\underline g}_f, g_\infty)$ with
${\underline g}_f\in G({\Bbb A}_f)$ and $g_\infty\in G(K_\infty)$. Elements of $G(K)$, viewed as
diagonally embedded in $G({\Bbb A})$, are not underlined.

\proclaim{(0.5) Definition} An automorphic form with respect to an open compact subgroup
$\frak K$ of $G(\Cal O)$ is a (complex-valued) function 
$f\: G({\Bbb A}) \to \Bbb C$ such that, for all\linebreak $\gamma \in G(K)$,
${\underline g}\in G({\Bbb A})$ and ${\underline k}\in {\frak K}Z(K_\infty)$, the equality
$f(\gamma {\underline g}{\underline k})=f({\underline g})$ holds.
 
 Moreover, it is called a cusp form if for all ${\underline g}\in G({\Bbb A})$,
$$\int_{K\backslash \Bbb A}f\bigg (\pmatrix 1&{\underline u}\\ 0&1\endpmatrix {\underline g}
\bigg )  d{\underline u}= 0$$
($d{\underline u}$ is the normalized Haar measure on the compact group $K\backslash \Bbb A$).
\endproclaim

  These notions were first used intensively in positive characteristic by V.G. Drinfeld
(\cite {Dr}), although many of its main properties for general reductive groups, were given
by G. Harder (\cite {Ha}). Recall that cusp forms that transform like the special representation
led to a Galois reciprocity law (Drinfeld, loc. cit.; see also \cite {vdP-Re}).\par
  It is clear that the field $\Bbb C$ does not play any role in this definition of an automorphic
form, so one may replace it by any commutative ring of characteristic zero (with unit). 
Moreover, the integral in
this definition is indeed a finite sum. So in the definition of a cusp
form, one can replace $\Bbb C$ by any commutative ring containing $\Bbb Q$ 
(as subring with unit).\par
  
\noindent {\bf (0.6)} In what follows,
 the compact subgroups $\frak K$ of $G(\Cal O)$ will be of the
form ${\frak K}={\frak K_f}\times {\frak K_\infty}$, with ${\frak K_f}$ and ${\frak K_\infty}$
 open compact subgroups of $G({\Cal O}_f)$ and $G({\Cal O}_\infty)$ respectively.

\vskip5pt Modular forms also exist in positive characteristic. They were introduced concretely by
E.-U. Gekeler and D. Goss (\cite {Gek1}, \cite {Go}, see also \cite {Co}). 
Their definition, close to the classical one, will be given in (1.8). 
For our purposes, we just recall that they are functions defined on the Drinfeld upper half-plane 
$\Omega ={\Bbb P}_C^1(C)-{\Bbb P}_C^1(K_\infty )$ 
with values in $C$ ($C$ is the completion of an algebraic closure of $K_\infty$). 

\vskip5pt It is possible to make a parallel with the classical case: $K$ with $\infty$ is the
analog of ${\Bbb Q}$ equipped with its ordinary absolute value, $K_\infty $ and $C$ are the
analogs of $\Bbb R$ and $\Bbb C$ respectively, $A$ looks like $\Bbb Z$ and
$\Omega $ is the analog of the Poincar\'e half-plane.
This parallel ends here (in our context). There is no direct link 
between modular forms and automorphic forms in positive characteristic. On the contrary, it is
well known that the two notions of modular forms and automorphic forms coincide in the
classical case (see \cite {Gel}, \S 3). The main reason for the difference,
between modular forms and automorphic forms in positive characteristic, is probably that there
are no tools to go from $\Omega $ to $G({\Bbb A})$, then to translate functions defined on
$\Omega$ to functions defined on $G({\Bbb A})$ (this can be easily accomplished in the
classical case).\par
  We do not know how to pass from $\Omega $ to $G({\Bbb A})$, but they are both related
to the Bruhat-Tits tree $\tau $ of $G(K_\infty)$. On the one hand,  $\tau $ is isomorphic to the
intersection graph  of the analytic reduction of $\Omega$, viewed as a rigid analytic space over
$C$  ( \cite {Fe-vdP} ch.{\uppercase\expandafter{\romannumeral 5}}, 
\cite{Gek-Re} \S 1, see also \cite {vdP}). On the other hand , if
${\frak K}$ is  an open compact subgroup of  $G({\Bbb A})$ of the form 
${\frak K}={\frak K_f}\times {\frak K_\infty}$ (see (0.6)), where ${\frak K_\infty}$
is the stabilizer in $ G(K_\infty)$ of an edge of $\tau $, we have a one-to-one map between
$G(K)\backslash G({\Bbb A})/{\frak K}Z(K_\infty)$ (compare with (0.5))
and a finite disjoint union of quotients 
 of the set of edges of $\tau $ by arithmetic subgroups of $G(K)$ 
(this will be stated precisely in the next paragraph).\par
  There exist functions on the set of edges of $\tau $ that are of  particular interest, 
namely, the
harmonic cocycles. They were first introduced in our context by V.G. Drinfeld (\cite {Dr}),
who proved that, when they take {\it values in a field of characteristic zero}, they are indeed
the automorphic forms that transform like the special representation (this result appears in the
proof of his reciprocity law, loc. cit.; see also \cite {vdP-Re} and (1.13) below).\par
  Harmonic cocycles, more precisely a generalization of the above ones, were compared with
modular forms by P. Schneider in the $p$-adic context (\cite {Sc}) and by J. Teitelbaum in
positive characteristic.
In \cite {Te} (see also (1.9) below), J. Teitelbaum proves that spaces of harmonic cocycles taking
{\it values in characteristic} $p$ are isomorphic to the spaces of modular forms. It seems to be
difficult to lift directly these harmonic cocycles to characteristic zero and then,
using Drinfeld's result, to compare it with automorphic forms.\par
The first result comparing modular and automorphic
forms has appeared in \cite {Gek-Re}, \S (6.5) (recalled in (1.10)): it relies modular forms of
{\it weight 2}, doubly cuspidal, with cusp forms (using Teitelbaum's result, loc. cit.).\par

\vskip5pt {\it The purpose of this paper is to study the relationships between automorphic
forms and modular forms (in positive characteristic)}. Then, using Teitelbaum's result, we try to
interpret harmonic cocycles of equal characteristic (i.e. with values in characteristic $p$, the
same as the base field $K$) as automorphic forms.\par
  In \S 2 we introduce a notion of automorphic forms of equal characteristic, i.e. taking values
in spaces of the same characteristic $p$ as the global field $K$. In \S 2, we also introduce a
notion of special representation (of equal characteristic), which is a variant of the usual
one. Then we compare harmonic cocycles and automorphic forms, both of equal characteristic
(theorem (2.4)); indeed, we prove that the harmonic cocycles of equal characteristic are also, in
some sense, automorphic forms that transform like the special representations(see (2.11)).\par
  The automorphic forms of equal characteristic that we introduce in \S 2 
are, as it can be easily seen, the reduction
modulo $p$ of ``automorphic forms" taking values in spaces of characteristic zero. These latter
forms are not exactly automorphic forms in the sense of Drinfeld, because, they do not
satisfy to conditions at $\infty$, but since we work with automorphic forms that transform like
the special representations, conditions at $\infty$ are not essential. We obtain a result (theorem
(3.7)), which interprets {\it modular forms of characteristic $p$ and of  weight $n+2$}
(or harmonic cocycles of equal characteristic $p$ and of the same weight)
as functions with values in characteristic zero. For the weight $2$, it completes a result of
\cite {Gek-Re} \S (6.5) (see (3.9)).

\vskip5pt \noindent {\bf (0.7)} For general notions of rigid analytic geometry, we refer to
\cite {B-G-R}, \cite {Ger-vdP} and \cite {Fe-vdP}. The Bruhat-Tits tree of $G(K_\infty)$ is
defined and extensively studied in \cite{Se}, ch. 2.  All that is needed, concerning the analytic
structure of the Drinfeld upper half-plane and its links with the Bruhat-Tits tree of
$G(K_\infty)$, is explained in \cite {Gek-Re}, \S 1. The underlying objects and tools that are used
here are Drinfeld modules and Drinfeld modular schemes: the details can be found in \cite {A-B}.

\head 1 Modular forms and harmonic cocycles.\endhead

  Let $\Omega ={\Bbb P}^1_C(C)-{\Bbb P}^1_C{(K_\infty)}$ be the Drinfeld upper half-plane.\par

\vskip5pt
\noindent {\bf (1.1)} Let $\pi $ be a uniformizing parameter of $K_\infty$, with $\pi \in K$. For
$n\in \Bbb Z$ we write $D_n$ for the subset of $z\in \Omega$ that satisfy 
$\vert \pi \vert^{n+1}\leq \vert z\vert \leq \vert \pi \vert^n$ and 
$\vert z-\rho \pi ^n\vert \geq \vert \pi \vert^n$, 
$\vert z-\rho \pi ^{n+1}\vert \geq \vert \pi \vert^{n+1}$
for all $\rho \in {\Bbb F}(\infty)^\star$, where ${\Bbb F}(\infty)\hookrightarrow K_\infty$ is 
isomorphic to the residue field of $K$ at $\infty$.\par
  For all $z\in K_\infty$ and $n\in {\Bbb Z}$ we set $D_{(n, z)}=z+D_n$. Let $I$ be the set of 
$(n, z)$ with $n\in \Bbb Z$ and $z$ belonging to a set of representatives of 
$K_\infty/\pi ^{n+1}{\Cal O}_\infty$. Then we have $\Omega =\cup_{i\in I}D_i$; more precisely,
$(D_i)_{i\in I}$ is a pure covering of $\Omega $. We denote the corresponding analytic reduction
by $R:\Omega \rightarrow \bar \Omega $ ; 
$\bar \Omega $ is a tree of ${\Bbb P}^1_{{\Bbb F}(\infty)}$, these ${\Bbb P}^1_{{\Bbb F}(\infty)}$
are its irreducible components, each of them meeting
$\sharp ({\Bbb F}(\infty))+1$ others in ordinary double points which 
are rational over ${\Bbb F}(\infty)$, and any two of them having at most one common point. 
We denote the intersection graph of $\bar \Omega $ by $T$. 
An edge $e$ of $T$ corresponds to the intersection of two irreducible
components of $\bar \Omega$, $C_1$ and $C_2$ say. Let ${\bar \Omega}_e $ be  the subset of 
$\bar \Omega $ equal to $C_1\cup C_2$ minus their intersection points with the other irreducible
components $C\not= C_1,C_2$. Then 
$(R^{-1}({\bar \Omega }_e))_{e}$ is the previous pure covering $(D_i)_{i\in I}$, where $e$ runs
over the set of non oriented edges of $T$.  
 
\vskip5pt \noindent {\bf (1.2)} Let $\tau $ be the Bruhat-Tits tree of $G(K_\infty)$. It is
canonically $G(K_\infty)$-isomorphic to $T$ (see \cite {Gek-Re}, \S 1). 
Now, the term {\it edge} means {\it oriented edge}. Let $e$ be an 
edge of $\tau $ or $T$, then $e(0)$, resp. $e(1)$, is its origin, resp. its end point; $-e$
is the edge with the origin and the end point interchanged. 

\vskip5pt \noindent {\bf (1.3)} Let $n\in \Bbb N$ and let $L$ be a ring, containing $K_\infty $
as a subring if $n\not= 0$. The ring $L$ is supposed to be commutative with unit, its subrings
are supposed to have the same unit {\it as all rings and subrings shall be in this paper}.
 We denote  
the subspace of $L[X, Y]$ (the polynomial ring in two variables) of homogeneous
polynomials of degree $n$ by $V_n(L)$. 
It is a free $L$-module of rank $n+1$. It is equipped with a 
$G(K_\infty)$-action, trivial for \linebreak $n=0$, denoted by $\rho_n:G(K_\infty)\to GL(V_n(L))$,
and defined in the following way for $n>0$: let 
$g \in G(K_\infty)$ be such that $g^{-1} =\pmatrix a&b\\ c&d\endpmatrix$ and let
$j$ be an integer, $0\leq j\leq n$, then $\rho_n(g) (X^jY^{n-j})=(aX+bY)^j(cX+dY)^{n-j}$.
We set $V_n(L)^\star =Hom_L(V_n(L), L)$. The following definition was given in \cite {Te}.

\proclaim{(1.4) Definition} Let $n\geq 2$ be an integer.
An $L$-harmonic cocycle of weight $n$ is a function\par
\noindent $f: edges(\tau )\rightarrow V_{n-2}(L)^\star$ such that:\par
\noindent (i) $f(-e)=-f(e)$ for all edges of $\tau $,\par
\noindent (ii) $\sum_{e(0)=v}f(e)=0$ for all vertices $v$ of $\tau $, where the sum is taken over
the edges with the origin equal to $v$.
\endproclaim

\noindent  {\bf (1.5)} Let $\underline{H}^{n}(L)$ be the set of $L$-harmonic cocycles of
weight $n$. It is an $L$-module, equipped with the following $G(K_\infty)$-action:
for all $f\in \underline{H}^{n}(L)$, $g\in G(K_\infty)$ and for all edges $e$ of $\tau $,
$g(f)(e)=\rho^\star_{n-2}(g)(f(g^{-1}e))$ (where $\rho^\star_{n-2}$ is the representation on
$V_{n-2}(L)^\star$ induces by $\rho_{n-2}$). 
If $\Gamma $ is a subgroup of $G(K_\infty)$, we denote 
the submodule of elements in $\underline{H}^{n}(L)$ fixed under the
$\Gamma $-action (coming from that of $G(K_\infty)$) by $\underline{H}^{n}(L)^\Gamma $.
$\underline{H}_!^{n}(L)^\Gamma $, resp. 
$\underline{H}_{!!}^{n}(L)^\Gamma $, are the submodules of elements in 
$\underline{H}^{n}(L)^\Gamma $ with finite supports modulo $\Gamma $, resp. which are zero on the
cusps of $\Gamma $. We do not explain this notion of cusp here because we do not use it except 
in the two recalls just below. 

\vskip5pt \noindent {\bf (1.6)} A subgroup $\Gamma $ of $G(K)$ is said to be {\it arithmetic}
if $\Gamma \cap G(A)$ (see (0.1)) is commensurable with both $\Gamma $ and $G(A)$. Let $\Gamma $
be such an arithmetic subgroup, then the quotient graph $\Gamma \backslash \tau $
is the union of a finite planar graph without ends, denoted $(\Gamma \backslash \tau )^\circ$, and
of  finitely many half-lines $({\Cal L}_i)_{1\leq i\leq c}$ (\cite {Se}, ch.2, th.9, p.143). This
half-lines are the cusps of $\Gamma $. Following \cite {Se} (Ch.2, lemme 6, p.142) and \cite {Te}
(prop.3), we have

\proclaim{(1.7) Proposition} Let $\Gamma $ be an arithmetic subgroup.
For any cusp $({\Cal L}_i)$ of $\Gamma $, let $e_i$ be its ``first edge" (i.e. its
edge with origin in $(\Gamma \backslash \tau )^\circ$). Then\par
\noindent (i) for all $n\geq 2$ and any $f\in \underline{H}_!^{n}(L)^\Gamma $, the support of $f$
modulo $\Gamma$ is included in \linebreak
$$edges((\Gamma \backslash \tau )^\circ)\cup \{ e_i\}_{1\leq i\leq c}$$
\noindent (ii) if $p$ does not divide zero in $L$ ($p=char(K)$), we have
$\underline{H}_{!!}^{2}(L)^\Gamma = \underline{H}_!^{2}(L)^\Gamma$,\par
\noindent (iii) for all $n\geq 2$, if $p$ is equal to zero in $L$, we have
$\underline{H}_!^{n}(L)^\Gamma = \underline{H}^{n}(L)^\Gamma$.
\endproclaim

  We will now introduce the notion of Drinfeld modular form. It was first studied in \cite
{Go} and \cite {Gek1}. To be short, we do not explain all their properties (as in (1.8) below),
but all can be found in \cite {Gek2}, \cite {Co} and \cite {Gek-Re}, \S 2.\par
\vskip5pt 
\noindent {\bf (1.8)} Let $\Gamma $ be an arithmetic subgroup of $G(K)$, and let $n\geq 2$ 
and $m\geq 0$ be integers. Recall that $C$ is the completion of an algebraic closure of
$K_\infty$. A Drinfeld modular form of weight $n$ and type $m$ with respect to $\Gamma $ is a
function $f:\Omega \to  C$ that satisfies\par
\noindent (i) for all $\gamma =\pmatrix a&b\\ c&d\endpmatrix \in \Gamma $ and for all 
$z\in \Omega $, $f(\gamma z)=(det \gamma )^{-m}(cz+d)^nf(z)$;\par
\noindent (ii) f is holomorphic on $\Omega $;\par
\noindent (iii) f is holomorphic at the cusps of $\Gamma $.\par
Moreover, we say that\par 
\noindent (iv) a modular form $f$ with respect to $\Gamma $ is cuspidal, resp.
$i$ times cuspidal, if it has a zero, resp. a zero of order at least $i$, in all cusps 
of $\Gamma $.\par
We denote the $C$-vector space of modular forms of weight $n$ and type $m$,
with respect to $\Gamma $  by $M_{n, m}(\Gamma )$, 
and the subspace of those which are $i$ times
cuspidal by $M_{n, m}^i(\Gamma )$. We also set $M_{n, m}^\star(\Gamma )=\cup_{i\geq 1}M_{n,
m}^i(\Gamma )$. These spaces are of finite dimension (\cite {Gek2}).\par
  As mentioned in the introduction, we have the following results

\proclaim{(1.9) Theorem} (\cite {Te}, th. 16)  Let $n\geq 2$ be an integer, then
the $C$-vector spaces  $M_{n, 0}^\star(\Gamma )$ and
$\underline{H}^{n}(C)^\Gamma$ are canonically isomorphic.\endproclaim

\proclaim{(1.10) Theorem} (\cite {Gek-Re}, \S (6.5)) Let $M_{2, 1}^2(\Gamma , {\Bbb F}_p)$ be the
subspace of elements $f$ of $M_{2, 1}^2(\Gamma )$ 
such that the residues of the holomorphic forms $f(z)dz$ are in \linebreak
${\Bbb F}_p={\Bbb Z}/p{\Bbb Z}$, then we have
$$\underline{H}_{!}^{2}({\Bbb Z})^\Gamma @>\text{reduction mod. p}>>
\underline{H}_{!!}^{2}({\Bbb F}_p)^\Gamma \simeq M_{2, 1}^2(\Gamma , {\Bbb F}_p)$$
the first map being surjective.\endproclaim

The proofs of these two results use the notion of residue for holomorphic differentials defined
on $\Omega $, which was introduced by M. van der Put in \cite {Fe-vdP}, ch.I. A holomorphic
form on $\Omega $ possesses a residue for each $D_i$, $i\in I$ (see (1.1)); with the aid of
the residue theorem (\cite {Fe-vdP}, ch.I, \S 3) and because of the isomorphism between the two
trees $T$ and $\tau $ (see (1.2)), it gives a harmonic cocycle...\par
 With the aid of (1.10) one can also prove 

\proclaim{(1.11) Theorem} (\cite {Gek-Re}, th. (6.5.3)) $M_{2, 1}^2(\Gamma )$ and
$\underline{H}_{!!}^{2}(C)^\Gamma$ are naturally isomorphic.
\endproclaim 

\vskip5pt 
\noindent {\bf (1.12)} A comparison theorem between automorphic forms and harmonic cocycles was
given by V.G. Drinfeld in the proof of his Galois reciprocity law (\cite {Dr}, see also
\cite {vdP-Re}, prop. 2.11). We now describe it.\par
  Let ${\frak K}_f$ be a an open compact subgroup of $G({\Cal O}_f)$. Then 
$G(K)\backslash G({\Bbb A}_f)/{\frak K}_f$ is finite; let $X\subset G({\Bbb A}_f)$ 
be a representative system of this set of double classes. For all ${\underline x}\in X$, set
$\Gamma_{\underline x}=G(K)\cap x{\frak K}_fx^{-1}$. It is an arithmetic subgroup of $G(K)$. 
Let $L$ be a ring containing $\Bbb Z$ (resp. $\Bbb Q$), and let 
${\Cal W}^{{\frak K}_f}(L)$ (resp. ${\Cal W}_\circ^{{\frak K}_f}(L)$)
be the space of automorphic forms (with values in $L$) with respect to
an open compact subgroup ${\frak K}$ of $G({\Cal O})$ of the form 
${\frak K}={\frak K}_f\times {\frak K}_\infty$ (resp. which moreover have finite supports in
$G(K)\backslash G({\Bbb A})/{\frak K}Z(K_\infty)$), 
where ${\frak K}_\infty$ is an open compact
subgroup of $G({\Cal O}_\infty)$ (see (0.5) and its comments).
Following G. Harder,  
${\Cal W}_\circ^{{\frak K}_f}(L)$ is the space of $L$-valued cuspidal automorphic forms
with respect to ${\frak K}_f$ (\cite {Ha}, (1.2.3)).\par
 Let $Sp_0(L)$ be the space of functions ${\Bbb P}_C^1(K_\infty)\to L$ 
that are locally constant in the rigid analytic sense, modulo constant
functions (a more general definition and details will be given in the next chapter). The group
$G(K_\infty)$ acts on $Sp_0(L)$: we denote this action by $sp_0$. For $f\in Sp_0(L)$ and 
$g\in G(K_\infty)$, $sp_0(g)f$ is the function $u\mapsto f(ug)$; $sp_0$ is the so called special
representation. 

\proclaim{(1.13) Theorem} (V.G. Drinfeld) One has the $L$-linear isomorphisms
$$\prod_{{\underline x}\in X}{\underline H}^{2}(L)^{\Gamma_{\underline x}}
=  Hom_{L[G(K_\infty)]}(Sp_0(L), {\Cal W}^{{\frak K}_f}(L))$$
$$\prod_{{\underline x}\in X}{\underline H}_!^{2}(L)^{\Gamma_{\underline x}}
=  Hom_{L[G(K_\infty)]}(Sp_0(L), {\Cal W}_\circ^{{\frak K}_f}(L))$$
\endproclaim

Following this theorem one says that harmonic cocycles, of weight $2$ and with values in
characteristic zero, are {\it automorphic forms that transform like the special representation}.
  
\vskip5pt \noindent {\bf (1.14)} 
We now summarize quickly what is known. Let ${\Cal R}_C$ be a local topological ring, having
$\Bbb Z$ equipped with the $p$-adic topology as topological subring and having $C$ as residue
field. It follows from (1.10-13) (and since we have spaces of finite dimension, \cite {Ha})

$$Hom_{G(K_\infty)}(Sp_0({\Cal R}_C), {\Cal W}_\circ^{{\frak K}_f}({\Cal R}_C))\simeq
\prod_{{\underline x}\in X}{\underline H}_!^{2}({\Cal R}_C)^{\Gamma_{\underline x}}
@>u>> \prod_{{\underline x}\in X}{\underline H}_{!!}^{2}(C)^{\Gamma_{\underline x}}
\simeq \prod_{{\underline x}\in X}M_{2, 1}^2(\Gamma_{\underline x})$$
the map $u$ being the reduction is surjective.

\head 2.  Automorphic forms of equal characteristic.\endhead

  The goal of this chapter is to give an analog of Drinfeld's theorem (1.13) for harmonic
cocycles of any weight, then for harmonic cocycles with values in characteristic $p$.
It will give an interpretation of modular forms (see (1.9)).

\proclaim{(2.1) Definition} Let $L$ be a ring of characteristic $p$. Let ${\frak K}_f$ be an open
compact subgroup of $G({\Cal O}_f)$. A $L$-valued automorphic form with respect to ${\frak K}_f$ is
a function $f:G(\Bbb A)@>>>L$ such that \par
\noindent (i) for all $\gamma \in G(K)$, ${\underline g}\in G({\Bbb A})$
and ${\underline k}_f\in {\frak K}_f$, the equality 
$f(\gamma {\underline g}{\underline k}_f)=f({\underline g})$ holds;\par
\noindent (ii) there exists an open compact subgroup ${\frak K}_\infty$
of $G({\Cal O}_\infty)$ such that the support
of $f$ is finite in 
$G(K)\backslash G({\Bbb A})/({\frak K}_f\times ({\frak K}_\infty Z(K_\infty)))$.\par
  We denote by ${\Cal W}_!^{{\frak K}_f}(L)$ the space of these automorphic forms.
\endproclaim 

  We have choosen to require no condition at $\infty$, we will see later that indeed the contrary
is also possible (see (2.11)).

\proclaim{(2.2) Definition} Let $n\in \Bbb N$ and let $L$ be a ring of characteristic $p$
containing $K_\infty$ if $n>0$. Let ${\Cal F}_n(L)$ be the space of locally constant functions
\linebreak
${\Bbb P}_C^1(K_\infty) @>>> V_n(L)$ and denote by $Sp_n(L)$ its quotient by the set of constant
functions. The group $G(K_\infty)$ acts on $Sp_n(L)$, we denote by $sp_n$ this action: for all 
$h\in Sp_n(L)$ and $g\in G(K_\infty)$  one has
$sp_n(g)h:z\mapsto \rho_n(g)h(zg)$ (see (1.3)). We call
$sp_n$ the ($L$-valued) special representation of rank $n$.
\endproclaim

\noindent {\bf (2.3)} 
In this definition, ``locally finite" means that , for all $h\in Sp_n(L)$, there exists a finite
open covering $(U_i)_{1\leq i\leq r}$ of ${\Bbb P}_C^1(K_\infty)$ such that $h$ is constant on
each $U_i$. ${\Bbb P}_C^1(K_\infty)$ can be viewed as the set of ends of $\tau $, i.e. as the
set of equivalent classes of half-lines of $\tau $, two half-lines being equivalent if their
intersection contains infinitely many edges (see \cite {Se} ch.2, p.100-101). For an
(oriented) edge $e$ of $\tau $ denote by $U(e)$ the set of equivalent classes of half-lines
containing $e$, then $U(e)_{e\in edges(\tau )}$ is a basis of open subsets for the topology of
${\Bbb P}_C^1(K_\infty)$ and, for all function $f:{\Bbb P}_C^1(K_\infty)\to V_n(L)$,
locally constant, there exists edges $e_1$,...,$e_r$ of $\tau $ and $\lambda_1$,...,
$\lambda_r$ in $V_n(L)$ such that $f=\sum_{1\leq i\leq r}\lambda_i1_{U(e_i)}$
($1_{U(e_i)}$ is the characteristic function of $U(e_i)$).\par
  Note that we have a $G(K_\infty)$-isomorphism: $Sp_n(L)\simeq Sp_0(L)\otimes_LV_n(L)$.

\proclaim{(2.4) Theorem}  Let $n\in \Bbb N$ and let $L$ be a ring of characteristic $p$
containing $K_\infty$ if $n>0$. Let ${\frak K}_f$ be an open compact subgroup of ${\frak K}_f$
and $X\subset G({\Bbb A}_f)$ be a set of representatives of
$G(K)\backslash G({\Bbb A}_f)/{\frak K}_f$. For all ${\underline x}\in X$ set
$\Gamma_{\underline x}=G(K)\cap {\underline x}{\frak K}_f{\underline x}^{-1}$. 
Then we have a $L$-isomorphism 
$$\prod_{{\underline x}\in X}{\underline H}^{n+2}(L)^{\Gamma_{\underline x}}
\simeq  Hom_{L[G(K_\infty)]}(Sp_n(L), {\Cal W}_!^{{\frak K}_f}(L))$$
($G(K_\infty)$ acts on ${\Cal W}_!^{{\frak K}_f}(L)$ via the regular representation of 
$G({\Bbb A)}$).
\endproclaim

 The proof needs many steps.

\vskip5pt \noindent {\bf (2.5)} Let $E=E(G({\Bbb A}_f)/{\frak K}_f, L)$ be the set of functions 
$f:G({\Bbb A}_f)\to L$ right invariant under ${\frak K}_f$. An element of
${\underline H}^{n+2}(E)$ can be viewed as a function 
$\varphi :edges(\tau )\times G({\Bbb A}_f)\to V_n^\star(L)$, then one sees that
${\underline H}^{n+2}(E)$ is equipped with the following action of $G(K)$: for all 
$\gamma \in G(K)$, $e\in edges(\tau )$ and ${\underline g}\in G({\Bbb A}_f)$,
$\gamma (\varphi )(e,{\underline g})=\rho^\star_n(\gamma )
(\varphi (\gamma^{-1}e,\gamma^{-1}{\underline g}))$, where $\rho^\star$ is the action of
$G(K_\infty)$ on $V_n^\star(L)$ coming from that on $V_n(L)$ (see (1.3)).

\proclaim {(2.6) Lemma} One has a $L$-isomorphism
${\underline H}^{n+2}(E)^{G(K)}\simeq _L
\prod_{{\underline x}\in X}{\underline H}^{n+2}(L)^{\Gamma_{\underline x}}$. 
\endproclaim

\demo {Proof} Let $\varphi :edges(\tau )\times G({\Bbb A}_f)\to V_n^\star(L)$ be an element of 
${\underline H}^{n+2}(E)^{G(K)}$. One has 
$G({\Bbb A}_f)=\coprod_{{\underline x}\in X}G(K){\underline x}{\frak K}_f$ (disjoint union). Let
$$
\matrix {\underline H}^{n+2}(E)^{G(K)}&@>\Phi >>& 
\prod_{{\underline x}\in X}{\underline H}^{n+2}(L)^{\Gamma_{\underline x}}\\
\varphi &\mapsto &(\varphi _{\underline x})_{{\underline x}\in X}\endmatrix
$$
where $\varphi _{\underline x}=\varphi ( \ ,{\underline x})$. For an edge $e$ of $\tau $, for
${\underline x}\in X$ and $\gamma \in \Gamma_{\underline x} $ with
$\gamma ={\underline x}{\underline k}{\underline x}^{-1}$ , where ${\underline k}\in {\frak K}_f$
(see the definition of $\Gamma_{\underline x}$ in theorem (2.4)), one has
$\rho^\star_n(\gamma^{-1})\varphi_{\underline x}(\gamma e)=
\rho^\star_n(\gamma^{-1})\varphi (\gamma e,{\underline x})= \varphi (e,\gamma^{-1}{\underline x})$
because $\varphi $ is invariant under $G(K)$. Then
$\rho^\star_n(\gamma^{-1})\varphi_{\underline x}(\gamma e)=\varphi (e,{\underline x}{\underline k})
=\varphi (e,{\underline x})$ which proves that $\Phi $ is well defined. The inverse map is given
by $(\psi )_{\underline x}\mapsto ((e,\gamma {\underline x}_0{\underline k})\mapsto 
\rho^\star_n(\gamma )\psi_{{\underline x}_0}(\gamma^{-1}e))$.
\qed \enddemo

\proclaim {(2.7) Lemma} One has a $L$-isomorphism 
$${\underline H}^{n+2}(E)^{G(K)}\simeq _L Hom_L(Sp_n(L),E)^{G(K)}$$ 
the action of $G(K)$ on
$Hom_L(Sp_n(L),E)$ coming from that on $Sp_n(L)$ via $sp_n$ and on $E$.\endproclaim

\demo {Proof} One interprets elements of ${\underline H}^{n+2}(E)^{G(K)}$ as in (2.5).
An element $\zeta \in Hom_L(Sp_n(L),E)^{G(K)}$ can be viewed as a function
$\zeta :Sp_n(L)\times G({\Bbb A}_f)\to L$ and recall that the functions of the form 
$\lambda 1_{U(e)}$ for $\lambda \in V_n(L)$ generate $Sp_n(L)$ (see (2.3)), then one can define  
$$
\matrix {\underline H}^{n+2}(E)^{G(K)}&@>\Psi  >>&Hom_L(Sp_n(L),E)^{G(K)}\\
\varphi &\mapsto &
((\lambda 1_{U(e)},{\underline g})\mapsto \varphi (e, {\underline g})(\lambda )).\endmatrix
$$
The inverse map to $\zeta \in Hom_L(Sp_n(L),E)^{G(K)}$ assigns the function \linebreak
$edges(\tau )\times G({\Bbb A}_f)\to V_n^\star(L)$ which maps $(e,{\underline g})$ to
$\zeta (. \ 1_{U(e)},{\underline g})$. This is the expected isomorphism.
\qed \enddemo

\proclaim {(2.8) Lemma} Let ${\Cal W}_?^{{\frak K}_f}(L)$ be the set of functions satisfying the
assertions of definition (2.1) except {\it (ii)}. Then, one has an $L$-isomorphism
$$Hom_L(Sp_n(L),E)^{G(K)}\simeq_L Hom_{L[G(K_\infty)]}(Sp_n(L),{\Cal W}_?^{{\frak K}_f}(L)).$$
\endproclaim

\demo {Proof} To a function $\zeta \in Hom_L(Sp_n(L),E)^{G(K)}$, viewed as in the proof of (2.7),
one associates $\Theta (\zeta ):Sp_n(L)\to {\Cal W}_?^{{\frak K}_f}(L)$, such that for 
$f\in Sp_n(L)$, $\Theta (\zeta )(f)$ is the function $G({\Bbb A}_f)\times G(K_\infty)\to L$ which
maps $({\underline g}_f,{\underline g}_\infty)$ to 
$\zeta (sp_n({\underline g}_\infty)f)({\underline g}_f)$. It is easy to see that it gives the
desired isomorphism.
\qed \enddemo

\demo {{\bf (2.9)} End of the proof of theorem (2.4)} One has to prove that one can replace 
${\Cal W}_?^{{\frak K}_f}(L)$ by the space ${\Cal W}_!^{{\frak K}_f}(L)$ of definition (2.1).
One uses the notations of the proofs of the three previous lemmata. 
Let $(\varphi _{\underline x})_{{\underline x}\in X}\in \prod_{{\underline x}\in X}
{\underline H}^{n+2}(L)^{\Gamma_{\underline x}}$ and let $\Theta (\zeta )$ be its image in
$Hom_{G(K_\infty)}(Sp_n(L),{\Cal W}_?^{{\frak K}_f}(L))$ by the
composition of the three preceeding isomorphisms (see proof of (2.8)). Let $\lambda \in V_n(L)$
and set $w=\Theta (\zeta )(\lambda 1_{U(e)})$, then, $w$ is a map 
$G({\Bbb A}_f)\times G(K_\infty)\to L$. Choose $\gamma \in G(K)$, ${\underline x}\in X$,
${\underline k}_f\in {\frak K}_f$ and $g_\infty\in G(K_\infty)$. One has 
$w(\gamma {\underline x}{\underline k}_f,g_\infty)=
\rho_n^{\star}(\gamma )(\varphi _{\underline x}(\gamma^{-1}g_\infty e))
(\rho_n(g_\infty)\lambda)$. It follows that 
$w(\gamma {\underline x}{\underline k}_f,g_\infty)\not= 0$ implies 
$\gamma_{-1}g_\infty e\in supp(\varphi _{\underline x})$. There exists a finite set 
$S\subset G(K_\infty)$ such that $supp(\varphi _{\underline x})\subset 
\Gamma_{\underline x}S{\frak K}_\infty Z(K_\infty)$, where ${\frak K}_\infty$ is the stabilizer
of $e$ in $G({\Bbb O}_\infty)$ (and $Z$ is the center of $G$).Then
$$supp(w)\cap [(G(K){\underline x}{\frak K}_f)\times G(K_\infty)]
\subset (G(K){\underline x}{\frak K}_f)\times (S{\frak K}_\infty Z(K_\infty)).$$ 
It finishes the proof of (2.4).
\qed \enddemo

It follows from (1.9) and (2.4)

\proclaim {(2.10) Corollary} There exists a $C$-linear isomorphism
$$\prod_{{\underline x}\in X}M_{n+2,0}^{\star}(\Gamma_{\underline x})
\simeq  Hom_{C[G(K_\infty)]}(Sp_n(C), {\Cal W}_!^{{\frak K}_f}(C)).$$
\endproclaim

\vskip5pt \noindent {\bf (2.11) Remark.} It is possible to prescribe a condition at $\infty$ in the
definition (2.1); now we explain this. We continue (2.1) by adding the following condition:\par

\vskip3pt \noindent {\it (iii) Let S(f) be the $L$-submodule of
${\Cal W}_!^{{\frak K}_f}(L)$ generated by $f$ and ${\frak K}_\infty$ acting on $f$
via the regular representation (i.e. the action of $k_\infty \in {\frak K}_\infty$ on $f$ gives
the function on $G(\Bbb A)$ ${\underline g}\mapsto f({\underline g}k_\infty)$).  
Then, there should exist a $L$-submodule $Q(f)$ of
$Sp_n(L)$ and a $L$-morphism $\varepsilon (f):Q(f)\to S(f)$, such that: $Q(f)$ is 
stable under $sp_n({\frak K}_\infty)$ and is generated, as 
$L[sp_n({\frak K}_\infty)]$-module, by one element; $\varepsilon (f)$ is surjective
and ${\frak K}_\infty$-equivariant}.\par
Let ${\Cal W}_{!,\infty}^{{\frak K}_f}(L)$ be the
set of elements of ${\Cal W}_!^{{\frak K}_f}(L)$ which satify {\it (iii).} One has
$$Hom_{L[G(K_\infty)]}(Sp_n(L), {\Cal W}_!^{{\frak K}_f}(L))
\simeq_L Hom_{L[G(K_\infty)]}(Sp_n(L), {\Cal W}_{!,\infty}^{{\frak K}_f}(L)).$$
{\it Proof.} One continues with the notations of (2.4)-(2.9). Let again, as in (2.9),
$\Theta (\zeta )\in Hom_{G(K_\infty)}(Sp_n(C), {\Cal W}_!^{{\frak K}_f}(C))$ and $u\in Sp_n(L)$.
Let ${\underline g}_f\in G({\Bbb A}_f)$, $g_\infty\in G(K_\infty)$ and set 
$f=\Theta (\zeta )(u)$. One has $f(({\underline g}_f,g_\infty))=\zeta
(sp_n(g_\infty)(u))({\underline g}_f)$. There exists edges of $\tau $, $(e_i)_{1\leq i\leq r}$,
and elements of $V_n(L)$, $\lambda_{1\leq i\leq r}$, such that 
$u=\sum_{1\leq i\leq r}\lambda_ie_i$ and one can choose as subgroup ${\frak K}_\infty$ for $f$ the
intersection of the stabilizers in $G(K_\infty)$ of the $e_i$'s. Let 
$Q(f)=L[sp_n({\frak K}_\infty)]u$ , then one has ($k_\infty$ is in ${\frak K}_\infty$) 
$\varepsilon (f)(sp_n(k_\infty)u)=({\underline g}\mapsto f({\underline g}k_\infty))$.\qed 

\vskip3pt This property (2.11) permits to say that {\it harmonic cocycles of weight $n+2$ are
automorphic forms that transform like the special representation of rank $n$.}\par

\vskip5pt \noindent {\bf (2.12)} Let ${\Cal A}_!^{{\frak K}_f}(L)$ be the set of functions
$\psi :Sp_0(L)\times G({\Bbb A})\to L$ such that, for all $u\in Sp_0(L)$ and 
${\underline g}\in G({\Bbb A})$, $\psi (u,.)$ satisfies assertions {\it (i)} and {\it (ii)} of
(2.1) and $\psi (.,{\underline g})$ is $L$-linear. The group $G(K_\infty)$ acts on
${\Cal A}_!^{{\frak K}_f}(L)$: if $g_\infty \in G(K_\infty)$, one has 
$g_\infty (\psi )(u,{\underline g})=\psi (sp_0(g_\infty^{-1})u,{\underline g}g_\infty)$. The next
proposition gives a variant of theorem (2.4). 

\proclaim {(2.13) Proposition} One has a natural $L$-isomorphism
$$ Hom_{L[G(K_\infty])}(Sp_n(L), {\Cal W}_!^{{\frak K}_f}(L))\simeq_L
 Hom_{L[G(K_\infty)]}(V_n(L), {\Cal A}_!^{{\frak K}_f}(L)).$$
\endproclaim

\demo {Proof} Let $\varphi \in Hom_{L[G(K_\infty])}(Sp_n(L), {\Cal W}_!^{{\frak K}_f}(L))$, 
it can be viewed as a function
$\varphi_1:Sp_n(L)\times G({\Bbb A})\to L$, then (because 
$Sp_n(L)\simeq V_n(L)\otimes_L Sp_0(L)$ as $L[G(K_\infty)]$-modules) as a function
$\varphi_2:V_n(L)\times Sp_0(L)\times G({\Bbb A})\to L$ satisfying the following properties: for
all $v\in V_n(L)$, $u\in Sp_0(L)$, ${\underline g}\in G({\Bbb A})$ and
$g_\infty \in G(K_\infty)$\par
\noindent - $\varphi_2(v,u,.)$ satisfies {\it (i)} and {\it (ii)} of (2.1),\par
\noindent - $\varphi_2(.,.,{\underline g})$ is $L$-bilinear,\par
\noindent - $\varphi_2(\rho_n(g_\infty)v, sp_0(g_\infty)u, {\underline g})=
\varphi_2(v,u,{\underline g}g_\infty)$.\par
\noindent Clearly, the map $\varphi \mapsto ( v\mapsto \varphi_2(v,.,.))$ gives the desired
isomorphism.
\qed \enddemo

\head 3. Characteristic zero and characteristic $p>0$. \endhead

  In all this paragraph, $L$ is a {\it field of characteristic $p>0$}, $n\geq 0$ is an integer and
we suppose moreover that $K_\infty \subset L$ if $n>0$. We want to lift our harmonic cocycles (or
modular forms) to characteristic zero. One needs firstly

\proclaim {(3.1) Proposition} $V_n(L)$ is a cyclic $L[G(K_\infty)]$-module. \endproclaim

\demo {Proof} One can suppose that $n>0$.\par
\vskip3pt \noindent (3.2) 
Let $\Cal D$ be the set of integers $m$ and $mp^r-1$ with $0<m<p$ and $r>0$. Note that the
binomial coefficient   $\pmatrix n\\ i\endpmatrix$ is not zero modulo $p$ 
for all $i$, $0\leq i\leq n$, if and only if $n\in \Cal D$. 

\vskip3pt \noindent (3.3) Let $n>0$ be an integer and let 
$\alpha =max\{\beta \in {\Cal D} \ / \ \beta \leq n\}$.
It is easy to see that $\alpha \geq n/2$.

\vskip3pt Let $n$ and $\alpha $ be as before. For $a$ in $K_\infty^\star$ let $\gamma_a$ and 
$\delta_a$ be the two matrices such that $\gamma_a^{-1}=\pmatrix a&1\\ 1&0\endpmatrix$ and
$\delta_a^{-1}=\pmatrix a&1\\ 0&1\endpmatrix$. One has (see (1.3))
$$\alignat3
\gamma_a(X^\alpha Y^{n-\alpha })&=
\sum_{0\leq i\leq \alpha } \pmatrix \alpha \\ i\endpmatrix a^iX^{n-\alpha +i}Y^{\alpha -i}
&\in L[G(K_\infty)]X^\alpha Y^{n-\alpha }&&\\
\delta_a(X^\alpha Y^{n-\alpha })&=
\sum_{0\leq i\leq \alpha } \pmatrix \alpha \\ i\endpmatrix a^iX^iY^{n-i}
&\in L[G(K_\infty)]X^\alpha Y^{n-\alpha }&&
\endalignat
$$
for all $a$ in $K_\infty$.
As $\pmatrix \alpha \\ i\endpmatrix\not= 0$, it follows from the first formula that\linebreak  
$X^{n-i}Y^i\in L[G(K_\infty)]X^\alpha Y^{n-\alpha }$ and from
the second formula that \linebreak $X^iY^{n-i}\in L[G(K_\infty)]X^\alpha Y^{n-\alpha }$,
for all $i$, $0\leq i\leq \alpha $. As $\alpha \geq n/2$ (see (3.3)), one has proved
$$V_n(L)=L[G(K_\infty)]X^\alpha Y^{n-\alpha }.$$
\qed \enddemo

\vskip5pt \noindent {\bf (3.4)} Let ${\Cal R}_L$ be a local ring of characteristic zero with
maximal ideal ${\frak M}_L$ and residue field ${\Cal R}_L/{\frak M}_L=L$. One denotes by $s$ the
canonical morphism ${\Cal R}_L\to {\Cal R}_L/{\frak M}_L=L$.\par

\vskip5pt \noindent {\bf (3.5)} Let
${\Cal A}_?^{{\frak K}_f}({\Cal R}_L)$ be the set of functions 
$f:Sp_0({\Cal R}_L)\times G({\Bbb A})\to {\Cal R}_L$ such that, for all $u\in Sp_0({\Cal R}_L)$
and 
${\underline g}\in G({\Bbb A})$, $f(u,.)$ satisfies assertions {\it (i)} of
(2.1) and $f(.,{\underline g})$ is ${\Cal R}_L$-linear (see definition of
${\Cal W}_?^{{\frak K}_f}(L)$ given in (2.8)).\par
\vskip5pt \noindent {\bf (3.6)} Let 
$\widetilde {Hom}_{L,G(K_\infty)}(V_n(L), {\Cal A}_?^{{\frak K}_f}({\Cal R}_L))$ be the space
of $G(K_\infty)$-linear maps \linebreak
$\psi :V_n(L)\to {\Cal A}_?^{{\frak K}_f}({\Cal R}_L)$ such
that, for any $v\in V_n(L)$, $s\circ \psi (v)$ is $L$-linear (see (3.4)).

\proclaim {(3.7) Theorem}  One has a natural surjective map
$$\widetilde {Hom}_{L,G(K_\infty)}(V_n(L), \Cal A_?^{{\frak K}_f}({\Cal R}_L))
\longrightarrow Hom_{L[G(K_\infty)]}(V_n(L), {\Cal A}_!^{{\frak K}_f}(L)).$$
\endproclaim
\demo {Proof}  
Let ${\Cal A}_?^{{\frak K}_f}(L)$ be the space of functions 
$h :Sp_0(L)\times G({\Bbb A})\to L$ such that, for all $u\in Sp_0(L)$ and 
${\underline g}\in G({\Bbb A})$, $h(u,.)$ satisfies assertion {\it (i)} of
(2.1) and $h(.,{\underline g})$ is $L$-linear. Let $f$ be in 
$ {\Cal A}_?^{{\frak K}_f}({\Cal R}_L))$ 
and ${\underline g}$ be in $G({\Bbb A})$, note that the ${\Cal R}_L$-linearity of
$f(.,{\underline g})$ implies that $s\circ f(.,{\underline g})$ is zero on ${\frak M}_LV_n(L)$,
and this last sentence is equivalent to say that $s\circ f(u,{\underline g})=0$ if $u$ takes
values in  ${\frak M}_L$ (because $u$ takes finitely many values).\par
  Let $f\in {\Cal A}_?^{{\frak K}_f}({\Cal R}_L)$ and $u\in Sp_0(L)$, we have
just seen that $s\circ f(u,.)$ makes sense; it defines a map
${\Cal A}_?^{{\frak K}_f}({\Cal R}_L)\to {\Cal A}_?^{{\frak K}_f}(L)$, which induces a
morphism 
$$\widetilde {Hom}_{L,G(K_\infty)}(V_n(L), {\Cal A}_?^{{\frak K}_f}({\Cal R}_L))
\longrightarrow Hom_{L[G(K_\infty)]}(V_n(L), {\Cal A}_?^{{\frak K}_f}(L)) .$$
With (3.1), one sees
that this map is surjective. Finally, as in (2.9) (see also (2.13)) one proves that
$$Hom_{L[G(K_\infty)]}(V_n(L), {\Cal A}_?^{{\frak K}_f}(L))\simeq 
Hom_{L[G(K_\infty)]}(V_n(L), {\Cal A}_!^{{\frak K}_f}(L)).$$
\qed \enddemo

\vskip5pt \noindent {\bf (3.8)}
Theorem (3.7), with (2.4) and (1.9), implies, when $L=C$, i.e. when $L$ is equal to the
completion $C$ of an algebraic closure of $K_\infty$,  that one has the following
diagram
$$
\align 
\widetilde {Hom}_{C,G(K_\infty)}(V_n(C), {\Cal A}_?^{{\frak K}_f}({\Cal R}_C)) 
  &\to Hom_{C[G(K_\infty)]}(Sp_n(C), {\Cal W}_!^{{\frak K}_f}(C))\\
  &\simeq \prod_{{\underline x}\in X}{\underline H}^{n+2}(C)^{\Gamma_{\underline x}}\\
  &\simeq \prod_{{\underline x}\in X}M_{n+2, 0}^\star({\Gamma_{\underline x}} )
\endalign
$$
the first map being surjective. Then, one sees that
modular forms in characteristic $p$, or harmonic cocycles in equal characteristic $p$,
of weight $n+2$, are indeed essentially objects comming from the characteristic zero. When 
$n=0$, one has a more precise result, which completes \cite {Gek-Re}, \S (6.5) (recalled in
(1.10)).

\proclaim {(3.9) Corollary} Let ${\Cal R}$ be a local ring of characteristic zero with residue
field
$L$ of characteristic $p>0$ and let $\Gamma $ be an arithmetic subgroup of $G(K)$. Then, one has a
natural surjective ${\Cal R}$-morphism
$${\underline H}^2({\Cal R})^\Gamma @>>> {\underline H}^2(L)^\Gamma .$$
\endproclaim

\demo {Proof} Let ${\Cal W}_?^{{\frak K}_f}({\Cal R})$ be the set of functions 
$G({\Bbb A})\to {\Cal R}$ such that for all $\gamma \in G(K)$, ${\underline g}\in G({\Bbb A})$ and
${\underline k}_f\in {\frak K}_f$ the equality 
$f(\gamma {\underline g}{\underline k}_f)=f({\underline g})$ holds (see (0.5), (1.12) and (2.8)).
It is easy to prove, as in (2.9) that
$$Hom_{{\Cal R}[G(K_\infty)]}(Sp_0({\Cal R}), {\Cal W}_?^{{\frak K}_f}({\Cal R}))\simeq_{\Cal R}
Hom_{{\Cal R}[G(K_\infty)]}(Sp_0({\Cal R}), {\Cal W}^{{\frak K}_f}({\Cal R}))$$
As in (2.13), one has also
$$Hom_{{\Cal R}[G(K_\infty)]}(Sp_0({\Cal R}), {\Cal W}_?^{{\frak K}_f}({\Cal R}))\simeq_{\Cal R}
Hom_{{\Cal R}[G(K_\infty)]}(V_ 0({\Cal R}), {\Cal A}_?^{{\frak K}_f}({\Cal R})).$$
Since  $V_ 0({\Cal R})={\Cal R}$ and $V_ 0(L)=L$, with trivial actions of $G(K_\infty)$, 
it is clear that there exists a surjective map
$$Hom_{{\Cal R}[G(K_\infty)]}(V_ 0({\Cal R}), \tilde {\Cal A}_?^{{\frak K}_f}({\Cal R})) @>>>
\widetilde {Hom}_{L,G(K_\infty)}(V_0(L), \tilde {\Cal A}_?^{{\frak K}_f}({\Cal R})).$$
This last map, the two previous isomorphisms, theorem (3.7) and (2.13) give 
$$Hom_{{\Cal R}[G(K_\infty)]}(Sp_0({\Cal R}), {\Cal W}^{{\frak K}_f}({\Cal R})) @>>>
Hom_{L[G(K_\infty)]}(Sp_0(L), {\Cal W}_!^{{\frak K}_f}(L))$$
which is surjective and, together with theorems (1.13) and (2.4), it gives the desired result.
\qed \enddemo

\head 4. Some comments. \endhead

Let $n$ and $l$ be two non-negative integers. As before, $C$ is the completion of an algebraic
closure of $K_\infty$.\par
 One can twist the representations $\rho_n$, that is, one can consider
$V_n(C)$ equipped with the action $g\mapsto det(g)^l\rho_n(g)\in GL(V_n(C))$ of $G(K_\infty)$ (see
(1.3)). One denotes by 
$V_{n,l}(C)$ the space $V_n(C)$ equipped with this last action, one denotes also by
${\underline H}^{n+2,l}(C)$ the harmonic cocycles with values in $V_{n,l}(C)^\star$ (see (1.4)).
The isomorphism  of (1.9) is proved for $l=0$ in \cite {Te}, but, with exactly the same arguments,
it can be extended to all $l\geq 0$. Then, one can prove, as in (3.8) 
$$
\align 
\widetilde {Hom}_{C,G(K_\infty)}(V_{n,l}(C), {\Cal A}_?^{{\frak K}_f}({\Cal R}_C)) 
  &\to Hom_{C[G(K_\infty)]}(Sp_n(C), {\Cal W}_!^{{\frak K}_f}(C))\\
  &\simeq \prod_{{\underline x}\in X}{\underline H}^{n+2,l}(C)^{\Gamma_{\underline x}}\\
  &\simeq \prod_{{\underline x}\in X}M_{n+2, l}^\star({\Gamma_{\underline x}} )
\endalign
$$ 
(and recall that there do not exist modular forms of weight one, see \cite {Co}, th. (6.9.1)).

\vskip5pt Let $\xi :K_\infty^\star \to C^\star$ be a character and let $V_{n,\xi }(C)$ be
$V_n(C)$ equipped with the action $g\mapsto \xi(det(g))\rho_n(g)$ of $G(K_\infty)$. Let 
${\underline n}=(n_1,\cdots ,n_r)\in {\Bbb N}^r$ and let ${\underline \xi }=(\xi_1,\cdots ,\xi_r)$
where the $\xi_i :K_\infty^\star \to C^\star$ are characters. Set
$$V_{{\underline n},{\underline \xi }}(C)=V_{n_1,\xi_1}(C)\otimes_C\cdots \otimes_C
V_{n_r,\xi_r}(C).$$
Harmonic cocycles with values in this space make sense, and properties closed to (2.4) or (3.7)
can be proved, but they have no interpretation by ``usual modular forms".\par
  Let ${\Bbb F}$ be a finite subfield  of $C$. The group $G({\Bbb F})$ acts on  
$V_{{\underline n},{\underline \xi }}(C)$ (for characters ${\Bbb F}^\star\to C^\star$ and by the
same law as $G(K_\infty)$). If $0\leq n_i\leq p-1$ ($p$ is the characteristic of our fields),
these representations are, up to isomorphisms, the irreducible representations of $G({\Bbb F})$
(\cite {Ba-Li}). One does not know what sort of representations of $G(K_\infty)$ are 
$V_{{\underline n},{\underline \xi }}(C)$. It is easy to see that, if $p$ divides $n$, the
representation $V_n(C)$ of $G(K_\infty)$ is not irreducible: 
$\oplus_{0\leq j \leq n/p}CX^{pj}Y^{n-pj}$ is a subrepresentation. May be, $V_n(C)$ is
an irreducible representation of $G(K_\infty)$ if and only if $0\leq n<p$ or $n=mp^r-1$ with
$0<m<p$ and $r>0$ (see (3.2))?

\enddocument